n-dimensional links, their components, and their band-sums

\documentclass[a4paper,10pt]{article} 
\begin{document}
\def\g{\gamma}    \def\ve{\varepsilon} \def\vp{\varphi}
\def\si{\roman{sin}} \def\co{\roman{cos}} \def\i{\rightarrow}
\def\e{\hookrightarrow} \def\l{\longrightarrow} \def\ttt{\longmapsto}
\def\pa{ pair of $n$-knots} \def\f{\noindent } \def\nl{\newline }
\def\np{\newpage } \def\x{\times } \def\p{\bf Proof of  }
\def\te{^t \hskip-1mm } 
\def\o{ordinary sense slice 1-link} 
\def\oo{ordinary sense slice 1-link }

\title{ $n$-dimensional links, their components, \\ and their band-sums }
\author{  Eiji Ogasa
\footnote{
This research was partially supported by Research Fellowships
of the Promotion of Science for Young Scientists.
\nl{\bf Keyword:}
$n$-links, $n$-knots, band-sums, ribbon $n$-knots, ribbon $n$-links, 
the Arf invariant of $n$-links (resp. $n$-knots), 
the signature of $n$-links (resp. $n$-knots), 
}\\
 ogasa@hep-th.phys.s.u-tokyo.ac.jp\\
High  Energy Physics Theory Groups\\
Department of Physics \\
University of Tokyo \\ 
Hongo, Tokyo 113, JAPAN\\}
\date{}
\maketitle

\baselineskip11pt {\bf Abstract. }  
\f We prove the following results (1) (2) (3)  
on relations between $n$-links and their components. 

(1) Let $L=(L_1, L_2)$ be a $(4k+1)$-link ($4k+1\geq 5$). 
Then we have 

\hskip1cm Arf $L$=Arf$L_1$+Arf$L_2$.

(2) Let $L=(L_1, L_2)$ be a $(4k+3)$-link ($4k+3\geq3$). 
Then we have 

\hskip1cm $\sigma L$=$\sigma L_1$+$\sigma L_2$.

(3) Let $n\geq1$. Then there is a nonribbon $n$-link $L=(L_1, L_2)$ such that 

\hskip5mm
$L_i$ is a trivial knot. 



\f 
We prove the following results (4) (5) (6) (7)
on band-sums of $n$-links.

(4) Let $L=(L_1, L_2)$ be a $(4k+1)$-link ($4k+1\geq 5$). 

Let $K$ be a band-sum of $L$.
Then we have 

\hskip1cm Arf $K$=Arf$L_1$+Arf$L_2$.

(5) Let $L=(L_1, L_2)$ be a $(4k+3)$-link ($4k+3\geq3$). 

Let $K$ be a band-sum of $L$.
Then we have 

\hskip1cm $\sigma K$=$\sigma L_1$+$\sigma L_2$.

The above (4)(5) imply the following (6). 

(6) Let $2m+1\geq3$. There is a set of three $(2m+1)$-knots 
$K_0$, $K_1$, $K_2$ with

the following property: 
$K_0$ is not any band-sum of any $n$-link $L=(L_1, L_2)$ 

such that $L_i$ is equivalent to $K_i$ ($i=1,2$).

(7) Let $n\geq1$. Then there is an $n$-link $L=(L_1, L_2)$ 
such that   
$L_i$ is 

\hskip5mm
a trivial knot ($i=1,2$) 
and  that a band-sum of $L$ is a nonribbon knot.



\f We prove a 1-dimensional version of (1). 

(8) Let $L=(L_1, L_2)$ be a proper 1-link.  Then 

\hskip1cm Arf $L$

\hskip1cm =Arf $L_1$+Arf $L_2$
+$\frac{1}{2}\{\beta^*(L)$+mod4 $\{\frac{1}{2}lk (L)\}\}$

\hskip1cm =Arf $L_1$+Arf $L_2$+mod2 $\{\lambda (L)\}$, 

where $\beta^*(L)$ is the Saito-Sato-Levine invariant 

and $\lambda(L)$ is the Kirk-Livingston invariant.

\np
\section{ Introduction and main results}

We work in the smooth category. 

An  {\it (oriented) (ordered) $m$-component n-(dimensional) link}
 is a smooth, oriented submanifold $L=\{L_1,...,L_m\}$ $\subset$ $S^{n+2}$, 
which is the ordered disjoint union of $m$ manifolds, each PL homeomorphic 
to the $n$-sphere. If $m=1$, then $L$ is called a {\it knot}.  

We say that n-links $L$ and $L'$ are {\it equivalent} 
if there exists an orientation preserving diffeomorphism 
$f:S^{n+2}$ $\rightarrow$ $S^{n+2}$ such that $f(L)$=$L'$  and 
$f\vert_{L}:L\rightarrow L'$ is 
an orientation and order preserving diffeomorphism.    
If $n$-knot $K$ bounds a $(n+1)$-ball $\subset S^{n+2}$, 
then $K$ is called {\it a trivial ($n$-)knot}.

We say that $m$-component $n$-dimensional links, $L$ and $L'$, 
are said to be {\it (link-)concordant} or {\it  (link-)cobordant} 
if there is a smooth oriented submanifold 
$\widetilde{C}$=\{$C_1$ ,...,$C_m$\} $\subset$ $S^{n+2}\times [0,1]$,  
which meets the boundary transversely in $\partial \widetilde{C}$, 
is PL homeomorphic to  $L \times[0,1]$ 
and meets 
$S^{n+2}\times \{0\}$ in $L$ 
 (resp. $S^{n+2}\times \{1\}$ in $L'$ ). 
An $n$-link $L$ is called a {\it slice } link 
if $L$ is cobordant to a trivial link.

We prove:

\vskip3mm\noindent{\bf  Theorem  1.1. } {\it 
(1) Let $4k+1\geq5$. 
Let $L=(L_1, L_2)$ be a $(4k+1)$-link. Then we have 
$$\mathrm{Arf} L=\mathrm{Arf} L_1+\mathrm{Arf} L_2.$$ 

(2) Let $4k+3\geq3$. 
Let $L=(L_1, L_2)$ be a $(4k+3)$-link. 
Then we have 
$$\sigma L=\sigma L_1+\sigma L_2.$$ 
}\vskip2mm

\f{\bf Note.} In \S2  we review the {Arf} invariant  and the signature. 

\vskip3mm
Furthermore we prove 1-dimensional version of Theorem 1.1.

\vskip3mm\noindent{\bf  Proposition 1.2. }{\it 
Let $L=(L_1, L_2)$ be a 1-link. 
Suppose that the $\mathrm{Arf}$ invariants of 2-component 1-links are defined, 
that is,  that the linking numbers are even. 

(1)  $\mathrm{Arf} L$=
$\mathrm{Arf} L_1+\mathrm{Arf}L_2+
\frac{1}{2}\{\beta^*(L)+\mathrm{mod4}\{\frac{1}{2}lk (L)\}\}$, 

where $\beta^*(L)$ is the Saito-Sato-Levine invariant. 

(2)  
$\mathrm{Arf}L=\mathrm{Arf} L_1+\mathrm{Arf} L_2$+mod2$\{\lambda (L)\}$, 

where $\lambda(L)$ is the Kirk-Livingston invariant.
}\vskip2mm 

\f{\bf Note.} 
The Saito-Sato-Levine invariant is defined 
in  \cite{Si} 
The Kirk-Livingston invariant is defined 
in \cite{KL}. 

\vskip2mm 
In order to continue to state our main results, we need some more definitions. 
An $n$-link $L= ( L_1,..., L_m ) $ is called a {\it ribbon } $n$-link 
if $L$ satisfies the following properties. 

(1) There is a self-transverse immersion $f:D^{n+1}_1\amalg...\amalg D^{n+1}_m\i S^{n+2}$ 
such that $f(\partial D^{n+1}_i)=L_i$. 

(2) The singular point set $C$  $(\subset S^{n+2}$) 
of $f$ consists of double points. 

$C$ is a disjoint union of $n$-discs $D^{n}_i (i=1,...,k)$. 

(3) Put $f^{-1}(D^{n}_j)=D^{n}_{jB}\amalg D^{n}_{jS}$. 
The $n$-disc $D^n_{jS}$ is trivially embedded in the interior 
Int $D^{n+1}_\alpha$ 
of a $(n+1)$-disc component $D^{n+1}_\alpha$.  
The circle $\partial D^n_{jB}$ is trivially embedded in 
the boundary $\partial D^{n+1}_\beta$ of 
an $(n+1)$-disc component $D^{n+1}_\beta$. 
The $n$-disc $D^n_{jB}$ is trivially embedded 
in the $(n+1)$-disc component $D^{n+1}_\beta$. 
(Note that there are two cases, $\alpha=\beta$ and $\alpha\neq\beta$.)

It is well-known that it is easy to prove that 
all ribbon $n$-links are slice. 

It is natural to consider the following. 

\vskip3mm
\noindent{\bf Problem A. }   
(1) Is there a nonribbon $n$-link $L=(L_1, L_2)$ 
such that $L_i$ is a ribbon knot ($i=1,2$)?

(2)
Is there a nonribbon $n$-link $L=(L_1, L_2)$ 
such that $L_i$ is a trivial knot ($i=1,2$)?
 \vskip3mm

The $n=1$ case holds because the Hopf link is an example. 
In \cite{O1}, 
the author gave the affirmative answer to the $n=2$ case. 
In this paper we give the affirmative answer to the $n\geq3$ case. 
We prove: 

\vskip3mm\noindent{\bf  Theorem  1.3.}{\it   
(1) Let $n\geq1$. Then there is a nonribbon $n$-link $L=(L_1, L_2)$ 
such that $L_i$ is a trivial knot. 

Furthermore we have the following. 

(2) Let $2m+1\geq1$.
Then there is a nonslice $(2m+1)$-link  $L=(L_1, L_2)$ 
such that $L_i$ is a trivial knot. 
(Note that $L$ is nonribbon since $L$ is nonslice.)

(3) Let $n\geq2$. 
Then there is a slice and nonribbon $n$-link $L=(L_1, L_2)$ 
such that $L_i$ is a trivial knot. 
}\vskip3mm

We need some more definitions. 
Let $L=(L_1, L_2)$ be an $n$-link. An $n$-knot $K$  is called 
a {\it band-sum } (of the components $L_1$ and $L_2$ ) of
the 2-link $L$ along a {\it band} $h$ if we have:

\f
(1) There is an $(n+1)$-dimensional 1-handle $h$, 
which is attached to $L$, embedded in $S^4$. 

\f
(2) 
There are a point $p_1\in L_1$ and a point $p_2\in L_2$. 
We attach $h$ to $L_1\amalg L_2$ along $p_1\amalg p_2$. 
 $h\cap (L_1\cup L_2)$ is the attach part of $h$.  
Then we obtain an $n$-knot from $L_1$ and $L_2$ 
by this surgery. The $n$-knot is $K$.

The set $(K_0, K_1, K_2)$ is called a {\it triple} of $n$-knots 
if $K_i$ is an $n$-knot. 
A triple  $(K_0, K_1, K_2)$ of $n$-knots is said to be {\it band-realizable} 
if there is an $n$-link $L=(L_1, L_2)$ 
such that $K_1$ (resp. $K_2$) is equivalent to $L_1$ (resp. $L_2$) 
and that $K_0$ is a band-sum of $L$.

 Note:  
Suppose that a triple  $(K_0, K_1, K_2)$ of $n$-knots is band-realizable. 
Then $[K_0]=[K_1]+[K_2]$, 
where [X] represents an element in the homotopy sphere group $\Theta_n$.  
See \cite{KM}  
for $\Theta_n$.  

It is natural to consider the following. 

\vskip3mm
\noindent{\bf Problem B. } 
Let $K_0, K_1, K_2$ be arbitrary $n$-knots. 
Then is the triple $( K_0, K_1, K_2 )$ of $n$-knots 
band-realizable?
\vskip3mm

By using the results in 
\cite{HK}\cite{K}\cite{KT}, 
we can prove: we have the affirmative answer to the $n=1$ case. 
By using 
\cite{Y}, 
we can prove: 
if $K_0, K_1, K_2$ are ribbon $n$-knots ($n\geq2$), 
we have the affirmative answer. 
In \cite{O1}, 
the author proved: there is a nonribbon 2-link $L=(L_1, L_2)$ 
such that $L_i$ is the trivial knot and that 
a band-sum of $L$ is a nonribbon knot.  
`Nonribbon case' of the $n$-dimensional version ($n\geq2$) 
is not solved completely. 
Thus, in this paper, 
we consider the following problems C, D, 
which are the special cases of Problem B.

\vskip3mm
\noindent{\bf Problem C. }
Is there a set of three $n$-knots $K_0, K_1, K_2$ 
such that the triple $(K_1, K_2, K_3)$ 
is not band-realizable?
\vskip3mm 

In this paper we give the affirmative answer when $n$ is odd and $n\geq3$. 
(Theorem 1.4, 1.5.) 

\vskip3mm
\noindent{\bf Problem D. } 
(1) Is there 
a set of one nonribbon $n$-knot $K_0$ and two ribbon $n$-knots $K_1, K_2$ 
such that the triple $(K_1, K_2, K_3)$ is band-realizable?

(2) Is there 
a set of one nonribbon $n$-knot $K_0$ and two trivial $n$-knots $K_1, K_2$ 
such that the triple $(K_1, K_2, K_3)$ is band-realizable?
\vskip3mm

In \cite{O1}, 
the author gave the affirmative answer to the $n=2$ case. 
In this paper we give the affirmative answer to the $n\geq3$ case. 
( Theorem 1.6. )  

\vskip3mm\noindent{\bf  Theorem 1.4.}  {\it 
Let $2m+1\geq3$. There is a set of three $(2m+1)$-knots $K_0$, $K_1$, $K_2$
such that the triple $(K_0, K_1, K_2)$  is not band-realizable. 
}\vskip3mm

Theorem 1.4 is deduced from Theorem 1.5. 

\vskip3mm\noindent{\bf  Theorem 1.5. } {\it 
(1) Let $4k+1\geq5$. 
Let $L=(L_1, L_2)$ be a $(4k+1)$-link. 
Let $K$ be a band-sum of $L$. 
Then we have 
$$\mathrm{Arf} K=\mathrm{Arf} L_1+\mathrm{Arf} L_2.$$

(2) Let $4k+3\geq3$. 
Let $L=(L_1, L_2)$ be a $(4k+3)$-link. 
Let $K$ be a band-sum of $L$. 
Then we have 
$$\sigma K=\sigma L_1+\sigma L_2.$$ 
}

\vskip3mm\noindent{\bf  Theorem 1.6.}{\it 
(1)  Let $n\geq1$. Let $T$ be a trivial $n$-knot. 
Then there is a nonribbon $n$-knot  $K$ 
such that the triple $(K, T, T)$ is band-realizable. 

Furthermore we have the following.

(2)  Let $2m+1\geq1$. Let $T$ be a trivial $(2m+1)$-knot. 
Then there is a nonslice $(2m+1)$-knot $K$ 
such that the triple $(K, T, T)$ is band-realizable. 
(Note that $K$ is nonribbon if $K$ is nonslice.)

(3) Let $2m+1\geq3$. Let $T$ be a trivial $(2m+1)$-knot $T$.  
Then there is a slice and nonribbon $(2m+1)$-knot $K$ 
such that the triple $(K, T, T)$ is 
band-realizable. 
}\vskip3mm

\f{\bf Note.} All even dimensional knots are slice. 
(\cite{Ke}.)    
\vskip3mm

\cite{O5} includes the announcement of this paper.
 
Our organization is as follows: 

\S2 Seifert matrices, the signature and 
the {Arf} invariant of $n$-knots (resp. $n$-links)  

\S3 Some properties of band-sums

\S4 Proof of Theorem 1.1.(1)

\S5 Proof of Theorem 1.1.(2)

\S6 Proof of Proposition 1.2

\S7 Proof of Theorem 1.5

\S8 Proof of Theorem 1.4

\S9 Proof of Theorem 1.3.(2)

\S10 Proof of Theorem  1.3.(3)

\S11 Proof of Theorem 1.3.(1)

\S12 Proof of Theorem 1.6.(2)

\S13 Proof of Theorem 1.6.(3)

\S14 Proof of Theorem 1.6.(1)

\S15 Open problems  

In \S9, we give an alternative proof of 
one of the main theorems of 
\cite{Kawauchi}   
and that of 
\cite{CappellShaneson}.   
In \S10, we give a short proof of the main theorem of 
\cite{Hitt}. 

\section{
Seifert matrices, the signature and the {Arf} invariant 
of $n$-knots (resp. $n$-links) 
} 

See the $n=1$ case \cite{Kf}\cite{Murasugi}\cite{Rolfsen}. 
See the $n\geq2$ case \cite{L1} \cite{L2}. 

Let $K$ be a $(2m+1)$-knot $(2m+1\geq1)$. 
Let $V$ be a connected Seifert hypersurface of $K$. 
Note the orientation of $V$ is compatible with that of $K$. 
Let $x_1,...,x_\mu$ be $(m+1)$-cycles in $V$ 
which are basis of $H_{m+1}(V;\bf{Z})$/Tor. 
Push $x_i$ to the positive direction of the normal bundle of $V$. 
Call it $x_i^{+}$.  
A {\it Seifert matrices} of $K$ associated with $V$ represented by basis 
$x_1,...,x_\mu$ 
is a matrix $A=(a_{ij})=(lk(x_i,x_i^{+}))$. 
Then we have: 
$A-(-1)^{m}\cdot ^t \hskip-1mm A$ represents the map 
 $\{H(V;\bf{Z})$/Tor\} $\x \{H(V;\bf{Z})$/Tor\} $\rightarrow\bf{Z}$,  
which is defined by the intersection product.

The {\it signature} $\sigma(K)$ of $K$ is the signature of 
the matrix $A+^t \hskip-1mm A$.  Therefore, we have:

\vskip3mm\noindent{\bf  Claim.}{\it 
If $2m+1=4k+3(\geq3)$, 
the signature of $K$ coincides with the signature of $\hat V$,  
where $\hat V$ is the closed oriented manifold 
which we obtain 
by attaching a $(4k+4)$-dimensional 0-handle to $\partial V$. 
}\vskip3mm

Let $K$ be a $(4k+1)$-knot $(4k+1\geq1)$. 
We regard naturally 
$(H_{2k+1}(V;\bf{Z})/\mathrm{Tor}) \otimes {\bf Z_2}$ 
as a subgroup of $H_{2k+1}(V;{\bf Z_2})$.  
Then we can take basis  $x_1,...,x_\nu,y_1,...,y_\nu$ 
of $(H_{2k+1}(V;\bf{Z})/\mathrm{Tor}) \otimes {\bf Z_2}$ 
such that 
$x_i\cdot x_j=0$, $y_i\cdot y_j=0$,  
$x_i\cdot y_j=\delta_{ij}$ 
for any pair $(i,j)$, 
where $\cdot$ is the intersection product. 
The {\it Arf invariant} of $K$ is mod 2 
$\Sigma_{i=1}^\nu lk(x_i,x_i^{+})\cdot lk(y_i,y_i^{+})$.


Let $L=(K_1, K_2)$ be a $(2m+1)$-link ($2m+1\geq1$). 
Let $V$ be a Seifert hypersurface of $L$. 
We define  $x_i,x_i^{+},  A, \sigma L$ in the same manner. 
If $2m+1=4k+3(\geq3)$, 
then $\sigma L$ is the signature of the closed oriented manifold 
$\hat{V}$, 
where $\hat V$ is the closed oriented manifold which we obtain 
by attaching $(4k+4)$-dimensional 0-handles to $\partial V$.

Let $L=(L_1, L_2)$ be a $(4k+1)$-link $(4k+1\geq1)$.  
We define the {Arf} invariant of 
 $L=(L_1, L_2)$ $(4k+1\geq1)$.  
There are two cases. 

(1) Let $4k+1\geq5$. 
The Arf invariant of $L$ is defined in the same manner
as the knot case. 

(2) Let $4k+1=1$. 
The {Arf} invariant of $L=(L_1, L_2)$ is defined only if 
the linking number lk$(L_1, L_2)$ of $L$ is even. 
Then we can take basis  $x_1,...,x_\nu,y_1,...,y_\nu, z$  
of $H_1(V;\bf{Z})$/Tor 
such that 
$x_i\cdot x_j=0, x_i\cdot y_j=\delta_{ij}, y_i\cdot y_j=0, 
x_i\cdot z=0, y_i\cdot z=0, z\cdot z=0$. 
The {\it {Arf} invariant} of $L$ is 
mod 2 $\Sigma_{i=1}^\nu lk(x_i,x_i^{+})\cdot lk(y_i,y_i^{+})$.
See e.g. Appendix of \cite{Kirby}.

\section{Some properties of band-sums} 

In our proof of main results we use the following properties of band-sums. 

\vskip3mm\noindent{\bf  Proposition 3.1.} {\it 
Let $L=(L_1, L_2)$ be an $n$-link. 
Let $K$ be a band-sum of $L$ along a band $h$. 

(1) $\mathrm{Arf} K=\mathrm{Arf} L$. $(n=4k+1\geq5.)$

(2) $\sigma K=\sigma L$. $(n=4k+3\geq3.)$

(3) Knot cobordism class of $K$ is independent of the choice of $h$. 
 $(n\geq2.)$

(4) The following two equivalent conditions hold. 
 $(n\geq1.)$ 

\hskip1cm(i)  If $L$ is slice, then $K$ is slice. 

\hskip1cm(ii) If $K$ is nonslice, $L$ is nonslice. 

(5)  The following two equivalent conditions hold. 
 $(n\geq1$).

\hskip1cm(i)  If $L$ is ribbon, then $K$ is ribbon. 

\hskip1cm(ii) If $K$ is nonribbon, $L$ is nonribbon. 
}\vskip3mm

\f{\bf Proof of (1)(2)(3).  }
We need a lemma.

\vskip3mm\noindent{\bf Lemma 3.2.} {\it 
There is a Seifert hypersurface $V$ for $L$ such that 
 $V\cap h$ is the attach part of $h$. 
}\vskip3mm

\f{\bf Proof of Lemma 3.2.} 
Let $h\x[-1,1]$ be a tubular neighborhood of $h\subset S^{n+2}$. 
Suppose $h\x[-1,1]\cap L$ is the attach part of $h$.  
Then we have 
$[L]=0\in H_n(\overline{S^{n+2}-(h\x[-1,1])};\bf{Z})$. 
By the following Claim 3.3, the above Lemma 3.2 holds. 
Claim 3.3 is proved by an elementary obstruction theory. 
(The author gave a proof in Appendix of \cite{O2}.)

\vskip3mm\noindent{\bf  Claim 3.3. }{\it 
Let $X$ be a compact oriented ($n+2$)-manifold with boundary.
Let $M$ be a closed oriented $n$-submanifold $\subset X$. 
We do not suppose that $M\cap X=\phi$ nor that $M\cap X\neq\phi$.  
Let $[M]=0\in H_n(X;\bf{Z})$. 
Then there is a compact oriented ($n+1$)-manifold $W$ 
such that $\partial W=M$. 
}\vskip3mm

Suppose that $n$ is odd and that $n\geq3$. 
By Lemma 3.2, 
a Seifert matrix of $L$ is a Seifert matrix of $K$. 
By \cite{L1}, \cite{L2}, 
Proposition 3.1.(1), (2), hold. 
Furthermore  Proposition 3.1.(3) holds when $n$ is odd and $n\geq3$.

Suppose that $n$ is even. By \cite{Ke}, all even dimensional knots are slice. 
Hence Proposition 3.1.(3) holds when $n$ is even.

\f{\bf Proof of (5).  }
Proposition 3.1.(5) 
holds by the definition of ribbon links.

\f{\bf Proof of (4).  }
If $L=\{L_1,...,L_m\}\subset S^{n+2}=\partial B^{n+3}$ is a slice 
$n$-link, then there is a disjoint union of embedded $(n+1)$-discs,  
$\widetilde{D}=\{D_1 ,...,D_m$\} $\subset B^{n+3}$, 
such that 
$\widetilde{D}$  meets the boundary transversely in $\partial \widetilde{C}$ 
and that $\partial D_i=L_i$. 
$\widetilde{D}$  is called a set of {\it slice discs} for $L$. 
If $L$ is a knot, 
$\widetilde{D}=D_1$ is called a {\it slice disc} for $L=L_1$.

We prove (i).   
Let $L=(L_1, L_2)$ be embedded in $S^{2m+3}=\partial B^{2m+4}=B^{2m+4}$. 
Let $D_1^{2m+2}\amalg D_2^{2m+2}\subset B^{2m+4}$ 
be a set of slice discs for $L$. 
Note that $D_1^{2m+2}\cap D_2^{2m+2}=\phi$. 
Then we can regard $h$ is a $(2m+2)$-dimensional 1-handle 
which is attached to $D_1^{2m+2}\amalg D_2^{2m+2}$. 
Put $D=h\cup D_1^{2m+2}\cup D_1^{2m+2}$. 
Then we can make a slice disc for $K$ from $D$.

\section{ Proof of Theorem 1.1.(1)} 

\vskip3mm\noindent{\bf  Theorem  1.1. } {\it 
(1) Let $4k+1\geq5$. Let $L=(L_1, L_2)$ be a $(4k+1)$-link. Then we have 
$$\mathrm{Arf} L=\mathrm{Arf} L_1+\mathrm{Arf} L_2.$$ 
}\vskip3mm

\f{\bf Proof. }
We prove:  

\vskip3mm\noindent{\bf Lemma 4.1.}{\it 
Let $K$ be a $(4k+1)$-knot 
 $\subset S^{4k+3}=\partial B^{4k+4}\subset B^{4k+4}$ ($4k+1\geq5$).  

Suppose that there is a compact ($4k+2$)-manifold $M$ 
which is embedded in $B^{4k+4}$ with the following properties. 

(1) $M\cap \partial B=\partial M=K$.

(2) $M$ intersects $\partial B^{4k+4}$ transversely.

(3) $H_i(M;{\bf Z})\cong 
H_i(\overline{\natural^\xi S^1\times S^{4k+1}-D^{4k+2}};{\bf Z})$   
 for each $i$, where 
 $\xi$ is 
 
 a nonnegative integer and 
 $\natural^0 S^1\times S^{4k+1}=S^{4k+2}$. 
 
Then we have  $\mathrm{Arf}(K)=0$. 
}\vskip3mm

\f{\bf Proof of Lemma 4.1.}
 We first prove: 
\vskip3mm\noindent{\bf Claim.}{\it 
In order to prove Lemma 4.1, it suffices to prove 
the case where $K$ is a simple knot. 
}\vskip3mm

\vskip3mm
\f{\bf Note.} See \cite{L2} 
for simple knots. 
Recall: If an $(2w+1)$-knot $K$ is a simple knot, then 
there is a Seifert hypersurface $V$ for $K$ with the following propositions. 
(1) $\pi_i(V)=0$ $i\leq w$. 
(2) There are embedded spheres in $V$ such that 
the set of the homology classes of the spheres 
is a set of generators of $H_{w+1}(V;\bf{Z})$.

\vskip3mm
\f{\bf Proof of Claim. }
Take a collar neighborhood of $S^{4k+3}=\partial B^{4k+4}\subset B^{4k+4}$.   
Call it $S^{4k+3}\x[0,1]$. 
Push $M\cap(S^{4k+3}\x[0,1])$ into the inside.  

By \cite{L2}, 
there is an embedding 
$f:S^{4k+1}\x[0,1]\e S^{4k+3}\x[0,1]$ 
with the following properties. 

(1) 
$f(S^{4k+1}\x\{1\})$ in $S^{4k+3}\x\{1\}$ is $K$. 

(2)     
$f(S^{4k+1}\x\{0\})$ in $S^{4k+3}\x\{0\}$ is a simple knot $K'$. 

Then {Arf}$K$={Arf}$K'$ and 
$M\cup f(S^{4k+1}\x[0,1])$ satisfies (1), (2), and (3) in Lemma 4.1.  
This completes the proof of the above Claim.

We prove Lemma 4.1 in the case where $K$ is a simple knot. 
There is a Seifert hypersurface $V$ for $K$ 
with the following properties: (1) $\pi_i V=0 (1\leq i\leq 2k)$. 
(2) There are embedded spheres in $V$ such that 
the set of the homology classes of the spheres 
is a set of generators of $H_{2k+1}(V;\bf{Z})$. 

Then we have: 

$$
H_i(V;{\bf Z})\cong
\left\{
\begin{array}{cl}
0  & \mbox{for $i\neq2k+1,0$} \\
{{\bf Z}^{2\mu}} & \mbox{for $i=2k+1$,}
\end{array}
\right.
$$

$H_{2k+1}(V;{\bf Z})\otimes{\bf Z_2}\cong H_{2k+1}(V;{\bf Z_2})$, 

$H_{2k+1}(V\cup M;{\bf Z})\otimes{\bf Z_2}\cong H_{2k+1}(V\cup M;{\bf Z_2})$, 

and 

$$
H_i(V\cup M;{\bf Z_2})\cong
\left\{
\begin{array}{cl}
0  & \mbox{for $i=2k+2$} \\
H_{2k+1}(V;{\bf Z_2})& \mbox{for $i=2k+1$} \\
0  & \mbox{for $i=2k$.} \\
\end{array}
\right.
$$

By Claim 3.3, there is a compact oriented $(4k+3)$-submanifold 
$W\subset B^{4k+4}$ 
such that $\partial W=V\cup M$.

Take the Meyer-Vietoris exact sequence:  
$$H_i(V\cup M;{\bf Z_2})\rightarrow H_i(W;{\bf Z_2})
\rightarrow  H_i(W, V\cup M;{\bf Z_2})$$. 

Consider the following part of the above sequence:  

$H_{2k+2}(V\cup M;{\bf Z_2})\rightarrow H_{2k+2}(W;{\bf Z_2})
\rightarrow  H_{2k+2}(W, V\cup M;{\bf Z_2})\rightarrow  $

$H_{2k+1}(V\cup M;{\bf Z_2})
{\to} H_{2k+1}(W;{\bf Z_2})
\rightarrow  H_{2k+1}(W, V\cup M;{\bf Z_2})\rightarrow  $

$H_{2k}(V\cup M;{\bf Z_2})$. 


Therefore we have 

$0\rightarrow H_{2k+2}(W;{\bf Z_2})
\rightarrow  H_{2k+2}(W, V\cup M;{\bf Z_2})\rightarrow  $

${\bf Z_2^{2\mu}} 
{\to} H_{2k+1}(W;{\bf Z_2})
\rightarrow  H_{2k+1}(W, V\cup M;{\bf Z_2})\rightarrow  $

0.

By using the Poincar\'e duality and the universal coefficient theorem, 
we have  
$ H_{2k+2}(W;{\bf Z_2})\cong  H_{2k+1}(W, V\cup M;{\bf Z_2})$ and 
$ H_{2k+1}(W;{\bf Z_2})\cong  H_{2k+2}(W, V\cup M;{\bf Z_2})$. 

Hence there is a set of basis $x_1,...,x_\mu,y_1,...,y_\mu$  
$\in H_{2k+1}(V;{\bf Z_2})$
with the following properties.

(1) $x_i\cdot x_j=0$, $y_i\cdot y_j=0$, $x_i\cdot y_j=\delta_{ij}$, 
where $\cdot$ denote the intersection product. 

(2) Let $f$ be the above map 
$H_{2k+1}(V\cup M;{\bf Z_2})\to H_{2k+1}(W;{\bf Z_2})$. 
Then $f(x_i)=0$

(3) $x_i$ is represented by an embedded $(2k+1)$-sphere in $V$.

We prove: 
\vskip3mm\noindent{\bf Lemma.}{\it 
If mod 2 lk($x_i^+, x_i$)=0 for each $i$, then $\mathrm{Arf}K$=0, 
where $x_i^+$ is one in \S2. 
}\vskip3mm

{\bf Proof.} Put 
$p: H_{2k+1}(V;{\bf Z})\rightarrow  H_{2k+1}(V;{\bf Z_2})$. 
There is a basis 
$\bar{x_1},...,\bar{x_\mu},\bar{y_1},...,\bar{y_\mu}$   
$\in H_{2k+1}(V;{\bf Z_2})$
with the following properties. 

(1) $\bar{x_i}\cdot \bar{x_j}=0$, $\bar{y_i}\cdot \bar{y_j}=0$, 
$\bar{x_i}\cdot \bar{y_j}=\delta_{ij}$, 
where $\cdot$ denote the intersection product. 
 
(2) $\bar{x_i}=p(x_i)$.   $\bar{y_i}=p(y_i)$. 

Then 
 lk($\bar{x_i}^+, \bar{x_i}$)$\equiv$lk($x_i^+, x_i$)   mod 2 
and 
 lk($\bar{y_i}^+, \bar{y_i}$)$\equiv$lk($y_i^+, y_i$)   mod 2.  

{Arf}$K$=
mod 2 {$\Sigma_{i=1}^\mu$  lk($\bar{x_i}^+, \bar{x_i}$)$\cdot$ lk($\bar{y_i}^+, \bar{y_i}$)}
=
mod 2 {$\Sigma_{i=1}^\mu$  lk(${x_i}^+, {x_i}$)$\cdot$ lk(${y_i}^+, {y_i}$)}
Hence the above Lemma holds.

Let $\alpha$ be a $\bf{Z_2}$-($2k+2$)-chain in $W$ which bounds $x_i$.  
Let $\beta$ be a $\bf{Z_2}$-($2k+2$)-chain in $S^{4k+3}$ which bounds $x_i$.  
Then $\gamma=$$\alpha\cup\beta$ is a $\bf{Z_2}$-($2k+2$)-cycle in $B^{4k+4}$. 
We prove: 

\vskip3mm\noindent{\bf Claim.}{\it 
The $\bf{Z_2}$-intersection product $\gamma\cdot\gamma$ in $B^{4k+4}$ is  
mod 2 lk($x_i^+, x_i$). 
}\vskip3mm

\f{\bf Proof. } 
Push off $\alpha$ to the positive direction of the normal bundle of $W$ in $X$. Call it $\alpha^+$.  
 Note $\alpha^+$ bounds $x_i^+$.
By considering the collar neighborhood $S^{4k+3}\x[0,1]$,  
we have that 
the $\bf{Z_2}$-intersection product $\gamma\cdot\gamma$ is 
the mod 2 number of the points $\alpha^{+}\cap\beta.$ 

It holds that  mod 2 lk($x^+, x$) is 
the mod 2 number of the points $\alpha^{+}\cap\beta.$
Hence  $\gamma\cdot\gamma=$ mod 2 lk$(x^+, x)$.


\vskip3mm\noindent{\bf Claim.}{\it 
The $\bf{Z_2}$-intersection product $\gamma\cdot\gamma$ 
in $B^{4k+4}$ is zero. 
}\vskip3mm

\f{\bf Proof. } 
$H_{2k+2}(B;{\bf Z_2})=0$. Hence $\gamma\cdot\gamma=0$. 

This completes the proof of Lemma 4.1. 

We go back to the proof of Theorem 1.1.(1).

In \cite{O5}, the author proved the following. 
 \cite{O4} includes the announcement. 
 
\vskip3mm\noindent{\bf Theorem. (See \cite{O4} \cite{O5}.)}
{\it 
Let $L_0=(L_{0a}, L_{0b})$ be a $(4k+1)$-link ($4k+1\geq5$). 
Then there is a boundary link $L_1=(L_{1a}, L_{1b})$ 
and a compact oriented submanifold 
$P\amalg Q \subset S^{4k+3}\times[0,1]$
with the following properties. 

(1)
$P=S^{4k+1}\times[0,1]$. 
Put $\partial P=P_0\amalg P_1$. 

$Q=\overline{(S^1\times S^{4k+1})-B^{4k+2}-B^{4k+2}}$.   
Put $\partial Q=Q_0\amalg Q_1$.

(2) $P$ (resp. $Q$) is transverse to $S^{4k+3}\times\{0,1\}$. 

(3) 
$(f(P_i), f(Q_i))$ in $(S^{4k+1}\times\{i\})$  is a link $L_i$
($i=0, 1$), 
where 
$(P\amalg Q)\cap (S^{4k+1}\times\{i\})$ is $(f(P_i), f(Q_i))$.  
}\vskip3mm

In order to prove Theorem 1.1.(1), 
it is suffices to prove that 
Arf $L_0$ = Arf$L_{0a}$+ Arf $L_{0b}$.  

Since $L_{0a}$ is cobordant to $L_{1a}$, we have 
  Arf $L_{0a}$=Arf $L_{1a}$. 
  
  Take  $L_{0b}\sharp L_{1b}$. 
  By using $Q$, we can make a manifold like $M$ in Lemma 4.1
  for $L_{0b}\sharp L_{1b}$. 
  By Lemma 4.1, we have Arf $L_{0b}$=Arf $L_{1b}$.

Since $L_1$ is a boundary link, 
there is a Seifert surface $V_{1a}$ for $L_{1a}$ 
( resp. $V_{1b}$ for $L_{1b}$ ) such that $V_{1a}\cap V_{1b}=\phi$. 
Let $K_1$ be a band-sum of $L_1$ 
along a band $h$ such that  
$h\cap \{V_{1a}\amalg V_{1b}\}$ is the attach part of $h$.
By considering  $V_{1a}$, $V_{1b}$, and $h$, 
we have Arf $K_1$=Arf $L_{1a}+$Arf $L_{1b}$.

Let $K_0$ be a band-sum of $L_0$. 
By Proposition 3.1.(1), Arf$K_0$=Arf $L_0$.

Take $L_0$ and $-L_1^\star$ in $S^{4k+3}$ such that 
$L_0$ is embedded in a ball $B^{4k+3}$ 
and that $-L_1^\star$ is embedded in $S^{4k+3}-B^{4k+3}$.  
Make $K_0$ in $B^{4k+3}$. 
Make $-K_1^\star$ in $S^{4k+3}-B^{4k+3}$.  
Take a connected-sum $K_0\sharp (-K_1^\star)$.  
By using $P\amalg Q$ and band-sums, 
we can make a manifold like $M$ in Lemma 4.1
for $K_0\sharp (-K_1^\star)$. 
By Lemma 4.1, we have Arf $K_0$=Arf $K_1$. 

Hence 
Arf $L_0$

=Arf $K_0$

=Arf $K_1$

=Arf $L_{1a}+$Arf $L_{1b}$

=Arf $L_{0a}+$Arf $L_{0b}$. 

This completes the proof of Theorem 1.1.(1).

\f{\bf Note. }  
(1) If  $bP_{4k+2}\cong\bf{Z_2}$, the proof of Theorem 3.5.(1) is easy. 
See \cite{KM} for $bP_{4k+2}$.  
Because: 
An arbitrary $n$-knot bounds a Seifert hypersurface. 
An arbitrary Seifert hypersurface is a compact oriented 
parallelizable manifold. 
Therefore $[K_0], [K_1], [K_2]\in bP_{4k+2}\subset \Theta_{4k+1}$. 
 If  $bP_{4k+2}\cong\bf{Z_2}$, 
then the {Arf} invariant of $K_i$ as a manifold coincides with 
the {Arf} invariant of $K_i$ as a knot. 

(2) There are integers $k$ such that 
$bP_{4k+2}\cong1$. 
See \cite{Browder} \cite{KM}.

\section{ Proof of Theorem 1.1.(2)} 

\vskip3mm\noindent{\bf Theorem  1.1. } {\it 
(2) Let $4k+3\geq3$. 
Let $L=(L_1, L_2)$ be a $(4k+3)$-link. 
Then we have 
$$\sigma L=\sigma L_1+\sigma L_2.$$ 
}\vskip3mm

\f{\bf Proof. }
Let $L=(L_+, L_-)$ be a 
($4k+3$)-link ($4k+3\geq3$) $\subset S^{4k+5}$.  
Let $V$ (resp. $V_+, V_-$) be a Seifert hypersurface of $L$ (resp. $L_+, L_-$). Take $S^{4k+5}\x[-1,1]$. 
Regard $L=(L_+, L_-)$ as in $S^{4k+5}\x\{0\}$. 

Take $L_+\x[0,1]$ in $S^{4k+5}\x[0,1]$ 
so that 
$L_+\x\{t\}$ is embedded in $S^{4k+5}\x\{t\}$ 
and that 
$L_+\x\{0\}$ coincides with $L_+$ in $L$ in $S^{4k+5}\x\{0\}$.

Take $L_-\x[-1,0]$ in $S^{4k+5}\x[-1,0]$ so that 
$L_-\x\{t\}$ is embedded in $S^{4k+5}\x\{t\}$ 
and that 
$L_-\x\{0\}$ coincides with $L_-$ in $L$ in $S^{4k+5}\x\{0\}$.

Then it holds that $(L_+\x\{0\}, L_-\x\{0\})$ in $S^{4k+5}\x\{0\}$  is $L$.

Take $V$ (resp. $V_+, V_-$) in 
$S^{4k+5}\x\{0\}$ (resp. $S^{4k+5}\x\{1\}$, $S^{4k+5}\x\{-1\}$).

Put 
$W=V_+\cup(L_+\x[0,1])\cup(-V)\cup (L_-\x[-1,0])\cup V_-$. 
Note 
$W\supset L$.

By Claim 4.2, 
there is a compact oriented ($4k+5$)-submanifold $X\subset S^{4k+5}\x[-1,1]$ 
such that $\partial X=W$. Hence 

$\sigma(W)=0$-----(i). 

By the definition of $W$,  

$\sigma(W)=$ $\sigma(V_+)+$$\sigma(-V)+$$\sigma(V_-)$
=$\sigma(V_+)-$$\sigma(V)+$$\sigma(V_-)$--------(ii). 

By (i)(ii),  $\sigma(V)=\sigma(V_+)+\sigma(V_-)$. 
Hence $\sigma(L)=\sigma(L_+)+\sigma(L_-)$.

\section{ Proof of Proposition 1.2}

\vskip3mm\noindent{\bf  Proposition 1.2. }{\it 
Let $L=(L_1, L_2)$ be a 1-link. 
Suppose that the $\mathrm{Arf}$ invariants of 2-component 1-links are defined, 
that is,  that the linking numbers are even. 

(1)  $\mathrm{Arf} L$=
$\mathrm{Arf} L_1+\mathrm{Arf}L_2+
\frac{1}{2}\{\beta^*(L)+\mathrm{mod4}\{\frac{1}{2}lk (L)\}\}$, 

where $\beta^*(L)$ is the Saito-Sato-Levine invariant. 

(2)  
$\mathrm{Arf}L=\mathrm{Arf} L_1+\mathrm{Arf} L_2$+mod2$\{\lambda (L)\}$, 

where $\lambda(L)$ is the Kirk-Livingston invariant.
}\vskip3mm

\f{\bf Proof. }
Put the Conway polynomial $\nabla_L(z)$ of $L=(L_1, L_2)$ to be 
$\nabla_L(z)=c_1\cdot z+c_3\cdot z^{3}+....$.
By Lemma 3.6 of \cite{Kf},  

\hskip3cm $c_1(L)=lk(L)$------------(i).

The Saito-Sato-Levine invariant $\beta(\quad)\in\bf{Z_4}$ 
 is defined in \cite{Si} 
for $L=(L_1, L_2)$ whose linking number is even.  
It is a generalization of the Sato-Levine invariant 
$\in\bf{Z_2}$ in 
\cite{Sato}.

Let lk($L$) be even.  
By Theorem 4.1 of 
\cite{Si}, 

$\beta^*(L)=$ $mod4 \{2c_3(L)-\frac{1}{2}c_1(L)\}$------------(ii).

By (i) and (ii), 

\hskip2cm
$\beta^*(L)$=mod4$\{2c_3(L)\}-$mod4$\{\frac{1}{2}lk(L)\}$--------(iii).

By (iii), 

\hskip2cm 
mod2$\{c_3(L)\}=$
$\frac{1}{2}\{\beta^*(L)+$mod4$\{\frac{1}{2}lk(L)\}\}$---(iv).

Note. The first $\frac{1}{2}$ in the right side make sense. 
We can regard the right side as an element in $\bf{Z_2}$.

The Kirk-Livingston invariant $\lambda(\quad)$ is defined in 
\cite{KL}.  
By the definition of $\lambda(\quad)$ and Theorem 6.3 of 
\cite{KL}, 
it holds that: If lk($L$) is even,

\hskip2cm mod4$\{\lambda(L)\}=$mod4$\{c_3(L)\}$.   ----------(v) 

By (v), 

\hskip2cm mod2$\{\lambda(L)\}=$mod2$\{c_3(L)\}$.   ----------(vi) 

By \cite{Murasugi}, 
it holds that: If lk($L$) is even, 

\hskip2cm mod2$\{c_3(L)\}$={Arf}$L$+{Arf}$L_1$+{Arf}$L_2$--------(vii).  


By (vi)(vii), Proposition 1.2.(2) holds.    
By (iv)(vii)Proposition 1.2.(1) holds.

\f{\bf Note. }
Let  lk($L$) be even. Then, by (iii) and (v),  we have 

\f$\beta^*(L)=$ mod 4$\{2\lambda(L)-\frac{1}{2}lk(L)\}$.  
It is written in  Addenda of [KL] that the author proved this result.

\section{Proof of Theorem 1.5} 

\vskip3mm\noindent{\bf  Theorem 1.5. } {\it 
(1) Let $4k+1\geq5$. 
Let $L=(L_1, L_2)$ be a $(4k+1)$-link. 
Let $K$ be a band-sum of $L$. 
Then we have 
$$\mathrm{Arf} K=\mathrm{Arf} L_1+\mathrm{Arf} L_2.$$

(2) Let $4k+3\geq3$. 
Let $L=(L_1, L_2)$ be a $(4k+3)$-link. 
Let $K$ be a band-sum of $L$. 
Then we have 
$$\sigma K=\sigma L_1+\sigma L_2.$$ 

}\vskip3mm

\f{\bf Proof of (1). }  
By Proposition 3.1(1), 
{Arf} $K$={Arf} $L$. 
 By Theorem 1.1.(1), 
$\mathrm{Arf} L=\mathrm{Arf} L_1+\mathrm{Arf} L_2$. 
Hence 
$\mathrm{Arf} K=\mathrm{Arf} L_1+\mathrm{Arf} L_2$.

\f{\bf Proof of (1). }  
By Proposition 3.1.(2), 
$\sigma K=\sigma L$. 
By Theorem 1.1.(2), 
$\sigma L=\sigma L_1+\sigma L_2$. 
Hence $\sigma K=\sigma L_1+\sigma L_2$.

\section{Proof of Theorem 1.4} 

\vskip3mm\noindent{\bf  Theorem 1.4.} {\it  
Let $2m+1\geq3$. There is a set of three $(2m+1)$-knots $K_0$, $K_1$, $K_2$
such that the triple $(K_0, K_1, K_2)$ is not band-realizable. 
}\vskip3mm

\f{\bf Proof of the ${\bf 2m+1=4k+1\geq5}$ case.}
There is a $(4k+1)$-knot ($4k+1\geq5$) whose {Arf} invariant is zero 
(resp. nonzero).

\f{\bf Proof of the ${\bf 2m+1=4k+3\geq3}$ case.} 
There is a $(4k+3)$-knot ($4k+3\geq3$) whose signature is zero (resp. nonzero).  


\section{Proof of Theorem 1.3.(2)} 

\vskip3mm\noindent{\bf  Theorem  1.3.}{\it   
(2) Let $2m+1\geq1$.
Then there is a nonslice $(2m+1)$-link  $L=(L_1, L_2)$ 
such that $L_i$ is a trivial knot. 
}\vskip3mm


\f{\bf Proof of the ${\bf 2m+1=4k+1(\geq1)}$ case. }
We prove:

\vskip3mm\noindent{\bf Proposition 9.1.}{\it 
There is a nonslice $(4k+1)$-link $L=(L_1,L_2)$ ($4k+1\geq 1 $) 
such that $L_i$ is a trivial knot ($i=1,2$).  
}\vskip3mm


\f{\bf Proof.}  Let $V_i$ be a Seifert surface for $L_i$. 
Let $V_i\cong \overline{(S^{2k+1}\x S^{2k+1})-B^{4k+2}}$. 
Suppose $V_1\cap V_2=\phi$. 
Let $a, b$ be basis of $H_{2k+1}(V_1, \bf{Z})$. 
Let $c, d$ be basis of $H_{2k+1}(V_2, \bf{Z})$. 

Let $K$ be a band-sum of $L$ along a band $h$. 
Suppose $h\cap V$ is the attach part of $h$.  
Put $V=V_1\cup V_2\cup h$.

We can suppose that a Seifert matrix of $L_1$ associated with $V_1$ 
represented by basis $a, b$ is

$\left(
\begin{array}{cc}
1&1\\
0&0\\
\end{array}
\right)$. 

We can suppose that a Seifert matrix of $L_2$ associated with $V_2$ 
represented by basis $c, d $ is 

$
\left(
\begin{array}{cc}
0&1\\
0&1\\
\end{array}
\right)
$

We can suppose that a Seifert matrix of $K$ associated with $V$ 
represented by basis $a, b, c, d $ is 

$\left(
\begin{array}{cccc}
 1&1&0&0\\
 0&0&1&0\\
 0&1&0&1\\
 0&0&0&1\\
\end{array}
\right)$.

One way of construction of $K$ is the following one: 
Take a ball $B^{4k+3}\subset S^{4k+3}$. 
Take a submanifold $V'_1$ in $B^{4k+3}$ which is equivalent to $V_1$. 
Take a submanifold $V'_2$ in $S^{4k+3}-B^{4k+3}$
which is equivalent to $V_2$. 
Let $L'_i$ be $\partial V'_i$. 
Take a connected-sum $L'_1\sharp L'_2$. 
By using pass-moves, we can make $K$ from $L'_1\sharp L'_2$. 
(Pass-moves for 1-knots are defined in 
\cite{Kf}. 
Pass-moves for $(2n+1)$-knots are defined by the author in 
\cite{O3}.($2n+1\geq3.$))  
 
we can make $K$ from $L'_1\sharp L'_2$.

We have 

det ($A+^{t}A$)=det
$
\left(
\begin{array}{cccc}
 2&1&0&0\\
 1&0&2&0\\
 0&2&0&1\\
 0&0&1&2\\
\end{array}
\right)
$
$=-15. $

$A+^{t}A$ is a ($4\x 4$)-matrix. 
Hence $\sigma(A+^{t}A)\neq 0$. 
Hence $K$ is nonslice.  
By Proposition 3.1.(4), $L$ is nonslice.  
This completes the proof when $2m+1=4k+1(\geq5)$.

\f{\bf Proof of the ${\bf 2m+1=4k+3(\geq3)}$ case.} 
We prove: 

\vskip3mm\noindent{\bf Proposition 9.2.}{\it 
There is a nonslice $(4k+3)$-link $L=(L_1, L_2)$ ($4k+3\geq3$) 
such that $L_i$ is a trivial knot ($i=1,2$).  
}\vskip3mm

\f{\bf Proof.}   Let $V_i$ be a Seifert surface for $L_i$. 
Let $V_i\cong \overline{(S^{2k+2}\x S^{2k+2})-B^{4k+4}}$. 
Suppose $V_1\cap V_2=\phi$. 

Let $a, b$ be basis of $H_{2k+2}(V_1, \bf{Z})$. 
Let $c, d$ be basis of $H_{2k+2}(V_2, \bf{Z})$.

Let $K$ be a band-sum of $L$ along a band $h$. 
Suppose $h\cap V$ is the attach part of $h$.  
Put $V=V_1\cup V_2\cup h$.

We can suppose that a Seifert matrix of $L_1$ associated with $V_1$ 
represented by basis $a, b$ is 

$
\left(
\begin{array}{cc}
1&1\\
0&0\\
\end{array}
\right)
$. 

We can suppose that a Seifert matrix of $L_2$ associated with $V_2$ 
represented by basis $c, d $ is 

$
\left(
\begin{array}{cc}
0&1\\
0&1\\
\end{array}
\right)
$.

We can suppose that a Seifert matrix of $K$ associated with $V$ 
represented by basis $a, b, c, d $ is

$A=
\left(
\begin{array}{cccc}
 1&1&0&0\\
 0&0&1&0\\
 0&-1&0&1\\
 0&0&0&1\\
\end{array}
\right)
$. 

We can construct $K$ 
by a similar way 
to the way of construction of the knot $K$ in Proof of Proposition 9.1.

By \S25, 26 of \cite{L2}, 
 $A$ is not a Seifert matrix of any slice knot. 
Hence $K$ is nonslice.  
By Proposition 3.1.(4), $L$ is nonslice.  
This completes the proof when $2m+1=4k+3(\geq3)$. 

This completes the proof of Theorem 1.3.(2).

\f{\bf Note. }
(1)  The above $(4k+3)$-knots $K$ are discussed in 
\cite{L2}\cite{Mio}. 

(2) By using the above link $L$, 
we can give a short alternative proof to one of the main results 
of 
\cite{CappellShaneson}, 
\cite{Kawauchi}. 
 The theorem is that there is a boundary 
 $(2m+1)$-link ($2m+1\geq1$) which is not cobordant to any split link.
Proof: If the above link $L$ is concordant to a split link, then $L$ is slice. 
Therefore $L$ is a boundary link which is not cobordant to any split link. 

\cite{Mio} 
prove a theorem which is close to this theorem 
but different from this theorem.

(3) We give a question: Do we give some answers to Problems in \S1 
by using 
\cite{CO}\cite{GL}\cite{L3}?   

\section{ Proof of Theorem 1.3.(3)}

\vskip3mm\noindent{\bf  Theorem  1.3.}{\it   
(3) Let $n\geq2$. Then there is a slice and nonribbon $n$-link 
$L=(L_1, L_2)$ such that $L_i$ is a trivial knot. 
}\vskip3mm

\f{\bf Proof of the ${\bf n\geq3}$ case.} 
Recall that the following facts hold by Theorem 4,1 of 
\cite{G} 
or by using the Mayer-Vietoris exact sequence.    
See, e.g., \S 14 of 
\cite{L2}  
 for the Alexander polynomials.   
See, e.g., p.160 of 
[Rolfsen]  
 and 
\cite{L1}  
for the Alexander invariant.   
Let $\widetilde {X}_{K}$ denote 
the canonical infinite cyclic covering 
of the complement of the knot $K$.

\vskip3mm\noindent{\bf Theorem 10.1. (known) }{\it   
Let $K$ be a simple $(2k+1)$-knot $(k\geq1)$.   
Let $\Delta_K(t)$ be the Alexander polynomial of $K$.   
Suppose  the $(k+1)$- Alexander invariant  
$H_{k+1}(\widetilde {X}_K;{\bf Q})\cong$ 
 $\{({\bf Q[t, t^{-1}]})$/ $\delta^1_K(t)\}$  
 $\oplus .  .  .   \oplus$
 $\{({\bf Q[t, t^{-1}]})$/ $\delta^p_K(t)\}$.   
 Then  
 $\Delta_K(t)=$  
 $a\cdot t^b\cdot\delta^1_K(t)$ $\cdot.  .  .  \cdot$ $\delta^p_K(t)$ 
  for a rational number $a$ and an integer   $b$  
and we can put $\Delta_K(1)=1.$
}

\noindent{\bf Theorem 10.2. (known) }{\it   
Let $K^{(n+1)}$ be the spun knot of $K^{(n)}$ ($n\geq1$).   
Let $H_{k}(\widetilde {X}_{K^{(n)}};{\bf Q})$ 
(resp. $H_{k}(\widetilde {X}_{K^{(n+1)}};{\bf Q})$  )
denote the $k$-Alexander invariant of $K^{(n)}$ (resp. $K^{(n+1)}$ ).   
Suppose that 
$K^{(n)}$ bounds a Seifert hypersurface $V$ such that 
$H_1(V;{\bf Z})\cong0$.
Then 
$H_2(\widetilde {X^{(n+1)}};{\bf Q})\cong$ 
$H_2(\widetilde {X^{(n)}};{\bf Q})$. 
}

\noindent{\bf Proposition 10.3. (known) }{\it   
Let $K^{(n+1)}$ be the spun knot of $K^{(n)}$ ($n\geq1$).   
If $K^{(n)}$ has a  simply connected Seifert hypersurface, 
then $K^{(n+1)}$ has a  simply connected 
Seifert hypersurface.  
}\vskip3mm

We prove: 
\vskip3mm\noindent{\bf Proposition 10.4.}{\it  
Let $K$ be a ribbon $n$-knot $\subset S^{n+2}$ ($n\geq1$). 
Then $H_2(\widetilde{X_K};\bf{Q})$ does not have 
 ${\bf Q[t, t^{-1}]}$-torsion.  
}\vskip3mm

\f{\bf Proof. }
Since $K$ is ribbon, 
there is a Seifert hypersurface $V$ which is diffeomorphic to 
$\overline{S^1\x S^n-D^{n+1}}$. 
It holds that 
$
H_i(V;{\bf Q})\cong
\left\{
\begin{array}{cl}
{\bf Q}  & \mbox{for $i=1,n$} \\
0  & \mbox{for $i\neq1, n$.} \\
\end{array}
\right.
$

Let $N(K)$ be a tubular neighborhood of $K$ in $S^{n+2}$. 
Put  $X=\overline{S^{n+2}-N(K)}$.
The submanifold $V\cap X$ is called $V$ again. 
Let $N(V)$ be a tubular neighborhood of $V$ in $X$. 
Put $Y=X-N(V)$.  
By using the Mayer-Vietoris exact sequence, it holds that 
$
H_i(Y;{\bf Q})\cong
\left\{
\begin{array}{cl}
{\bf Q}  & \mbox{for $i=1,n$} \\
0  & \mbox{for $i\neq1, n$.} \\
\end{array}
\right.
$

Let $p:\widetilde{X_K}\to X$ be the canonical projection map. 
Put $p^{-1}(N(V))=\amalg_{j=-\infty}^{\infty} V_j$.  
Put $p^{-1}(Y)=\amalg_{j=-\infty}^{\infty} Y'_j$.  
Suppose $\partial Y'_j\subset V_j\amalg V_{j+1}$.
Put $Y_j=V_j\cup Y'_j\cup V_{j+1}$. 
Then there is the Mayer-Vietoris exact sequence: 

$$
H_i(\amalg_{j=-\infty}^{\infty} V_j; {\bf Q})\to 
H_i(\amalg_{j=-\infty}^{\infty} Y_j; {\bf Q})\to 
H_i(\widetilde{X_K}; {\bf Q}). 
$$

Consider the following part:  
$
H_2(\amalg_{j=-\infty}^{\infty} Y_j; {\bf Q})\to 
H_2(X_K; {\bf Q})\to  
H_1(\amalg_{j=-\infty}^{\infty} V_j; {\bf Q}). 
$

Hence $0\to H_2(X_K; {\bf Q})\to \oplus^\mu{\bf Q[t, t^{-1}]}$ 
is exact, 
where $\mu$ is a nonnegative integer. 
Therefore Proposition 10.4 holds.    

Take a 3-link $L=(L_1, L_2)$ $\subset S^5$ 
in the proof of Theorem 1.3.(2).  

Suppose $L\subset \bf{R^5}\subset S^5$. 
Let $\alpha:{\bf R^4}\x{\bf R}\rightarrow
 {\bf R^4}\x{\bf R}$ be the map defined by 
$(x, y)\mapsto(x, -y)$.

Suppose that $L\subset{\bf R^4}\x\{y\vert y\geq0\}$, 
that 
$L_i\cap({\bf R^4}\x\{y\vert y=0\})$ is a 3-disc $D^3_i$.   
Note $D^3_1\cap D^3_2=\phi$. 
The link $\alpha(L)$ is called $-L^*=(-L_1^*, -L_2^*)$. 
The link

\f$(\overline{\{L_1\cup (-L_1^*)\}-D^3_1}, 
    \overline{\{K_2\cup (-K_2^*)\}-D^3_2} )$
is called $\widetilde L=(\widetilde L_1, \widetilde L_2)$.

We prove:

\noindent{\bf Claim.}{\it 
 $\widetilde L$ is slice and nonribbon. 
}

\f{\bf Proof. }Firstly we prove that $\widetilde L$ is slice. 
Take ${\bf R^4}\x{\bf R}\x\{z\vert z\geq0\}$. 
Regard 
\f${\bf R^4}\x{\bf R}$ as ${\bf R^4}\x{\bf R}\x\{z\vert z=0\}$. 
Put 

\f$F_\theta=$ ${\bf R^4}\x\{(y,z)\vert y=r\cdot\mathrm{cos}\theta,  
z=r\cdot\mathrm{sin}\theta, r\geq0, \theta:$ fix. $\}$, 
where $0\leq\theta<\pi$.  
Regard 
${\bf R^4}\x{\bf R}\x\{z\vert z\geq0\}$ 
as the rotating of $F_0$ around the axis  
${\bf R^4}\x\{0\}\x\{0\}$. 
When rotating  $F_0$, we rotate  $\overline{L-D^3_1-D^3_2}$ as well. 
The result of rotating $\overline{L-D^3_1-D^3_2}$ is 
a set of slice discs for $\widetilde L$.
Hence $\widetilde L$ is slice.

Secondly we prove that $\widetilde L$ is nonribbon.
Take Seifert hypersurfaces $V_1$ and $V_2$ in 
the proof of Theorem 1.3.(2) for $L$. 
Suppose $V_1, V_2\subset {\bf R^4}\x\{y\vert y=0\}$. 
 Suppose that 
$V_i\cap({\bf R^4}\x\{y\vert y=0\})$=$D^3_i$. 
Put $\widetilde {V_i}=V_i\cup\alpha(V_i)$.

Let $\widetilde{K}$ be a band-sum of 
$\widetilde{L}=(\widetilde{K_1}, \widetilde{K_2})$ along a band $h$.
Suppose $h$ $\cap\widetilde {V_i}$ is the attach part of $h$. 
Then a Seifert matrix of the 3-knot $\widetilde{K}$ is 

$P=
\left(
\begin{array}{cc}
X&0\\
 0&-X\\
\end{array}
\right)
$, 
where  
$X=
\left(
\begin{array}{cccc}
 1&1&0&0\\
 0&0&1&0\\
 0&-1&0&1\\
 0&0&0&1\\
\end{array}
\right)
$. 

Then det$(X+t\cdot ^t \hskip-1mm X)$ is 
the Alexander polynomial of $\widetilde{K}$ (see 
\cite{L1}\cite{L2} 
). 
By Theorem 10.1, 
$H_{m+1}(X_{\widetilde{K}};{\bf Q})$ has  ${\bf Q[t, t^{-1}]}$-torsion. 
By Proposition 10.4, $\widetilde{K}$ is nonribbon.   
By Proposition 3.1.(5), $\widetilde{L}$ is nonribbon. 

Next we prove the $n\geq3$ case. 
Let $L^{(3)}=(L_1^{(3)}, L_2^{(3)})$ be the above 3-link 
$\widetilde{L}=(\widetilde{L_1}, \widetilde{L_2})$.
For any band-sum $K^{(3)}$ of $L^{(3)}$, 
$H_2(X_{K^{(3)}})$ has  ${\bf Q[t, t^{-1}]}$-torsion.

Let  $L^{(n+1)}=(L_1^{(n+1)}, L_2^{(n+1)})$ be a spun link of 
$L^{(n)}=(L_1^{(n)}, L_2^{(n)})$ ($n\geq3$). 
We can take $h^{(n)}$ and $h^{(n+1)}$ so that 
the band-sum $K^{(n+1)}$ of $L^{(n+1)}$ 
along $h^{(n+1)}$ is a spun knot of $K^{(n)}$. 
( Put the core of the band in the axis of the rotation. )

By Theorem 10.2 and Proposition 10.3, 
 $H_2(\widetilde{X_{K^{(n)}}};{\bf Q})\cong 
  H_2(\widetilde{X_{K^{(3)}}};{\bf Q})$ ($n\geq3$).  
Hence $H_2(X_{K^{(n)}})$ has  ${\bf Q[t, t^{-1}]}$-torsion. 
By Proposition 10.4,  $K^{(n)}$ is nonribbon ($n\geq3$). 
By Proposition 3.1.(5),  $L^{(n)}$ is nonribbon. 

Since $L^{(n)}$ is a spun link, $L^{(n)}$ is slice ($n\geq4$). 
Hence $L^{(n)}$ is slice ($n\geq3$).

This completes the proof of the $n\geq3$ case.

\f{\bf Proof of the ${\bf n=2}$ case  } 
In \cite{O1} 
the author made a nonribbon 2-link as follows: 
Let $K$ be a 2-knot. 
Let $N(K)$ be a tubular neighborhood. 
We made a way to construct a 2-link $L^K=(L^K_1, L^K_2)$ in $N(K)$. 
We proved that there is a 2-knot $K'$ 
such that $L^{K'}$ is nonribbon. 

We prove: $L^{K'}$ is slice. 
Because: Let $K\subset S^4=\partial B^5=B^5$.
Take a slice disc $D^K\subset B^5$ for $K$. 
Take a tubular neighborhood $N(D^K)$ of $D^K$.  
Note $N(D^K)\cap S^4=N(K).$ 
Suppose that $K$ is a trivial knot 
    and that $D^K$ is embedded trivially in $B^5$ . 
Then we can make a set of slice discs 
$(D^K_1, D^K_2)$ for $(L^K_1, L^K_2)$ 
such that $(D^K_1, D^K_2)$ is embedded in  $N(D^K)$. 
Take a diffeomorphism map $f:N(D^K)\to N(D^{K'})$ 
such that $f(L^K_i)=L^{K'}_i$. 
The submanifold $f(D^{K}_i)$ is called $D^{K'}_i$.
Then $(D^{K'}_1, D^{K'}_2)$ is a set of slice discs for $L^{K'}$.

\f{\bf Note.} 
(1)  By using this section we can give a short alternative proof 
of the main theorem of 
\cite{Hitt}: 
there is a nonribbon and slice $n$-knot ($n\geq3$). 

( Nonribbon 2-knots and nonribbon 1-knots are known 
before 
\cite{Hitt}   
is written as 
\cite{Hitt} quoted. ) 

(2) In Proposition 10.4, furthermore, we can prove that 
  $H_i(\widetilde{X_K};\bf{Z})=0$ for $2\leq i\leq n-1$.

\section{ Proof of Theorem 1.3.(1)}

\vskip3mm\noindent{\bf  Theorem  1.3.}{\it   
(1) Let $n\geq1$. Then there is a nonribbon $n$-link $L=(L_1, L_2)$ 
such that $L_i$ is a trivial knot. 
}\vskip3mm

\f{\bf Proof. }
The $n=1$ case holds because the Hopf link is an example. 
The $n\geq2$ case follows from Theorem 1.3.(2), (3). 
This completes the proof.

\section{ Proof of Theorem 1.6.(2)}

\vskip3mm\noindent{\bf  Theorem 1.6.}{\it 
(2)  Let $2m+1\geq1$. Let $T$ be a trivial $(2m+1)$-knot. 
Then there is a nonslice $(2m+1)$-knot $K$ 
such that the triple $(K, T, T)$ is band-realizable. 
}\vskip3mm

\f{\bf Proof. }
 $K$ and $L=(L_1, L_2)$ in the proof of Theorem 1.3.(2) 
give examples.

\section{ Proof of Theorem 1.6.(3)}
\vskip3mm\noindent{\bf  Theorem 1.6.}{\it 
(3) Let $2m+1\geq3$. Let $T$ be a trivial $(2m+1)$-knot $T$.  
Then there is a slice and nonribbon $(2m+1)$-knot $K$ 
such that the triple $(K, T, T)$ is band-realizable. 
}\vskip3mm

\f{\bf Proof. }  
$K$ and $L=(L_1, L_2)$ in the proof of Theorem 1.3.(3) give examples.

\section{ Proof of Theorem 1.6.(1)}

\vskip3mm\noindent{\bf  Theorem 1.6.}{\it 
(1)  Let $n\geq1$. Let $T$ be a trivial $n$-knot. 
Then there is a nonribbon $n$-knot  $K$ 
such that the triple $(K, T, T)$ is band-realizable. 
}\vskip3mm

\f{\bf Proof. }
$K$ and $L=(L_1, L_2)$ in the proof of Theorem 1.3.(3) give examples
for the $n\geq3$ case. 
$K$ and $L=(L_1, L_2)$ in the proof of Theorem 1.3.(2) give examples
for the case where $n\geq1$ and $n$ is odd. 
The $n=2$ case follows from 
\cite{O1}. 
This completes the proof.

\section{ Open problems}
 
Even dimensional case of Problem C in \S1 is open. 
If the answer to the following problem is positive, 
then the $n=2$ case of Problem C is positive. 

\noindent{\bf Problem 15.1}
Let $L=(K_1, K_2)$ be a 2-link. 
Do we have: $\mu(L)=\mu(K_1)+\mu(K_2)$?  
  
See \cite{R}   
for the $\mu$-invariant of 2-knots. 
See \cite{O2}    %
for the $\mu$-invariant of 2-links. 

In \cite{O2} 
the author proved: 
if $L$ is a SHB link, the answer to Problem 15.1 is positive. 

If the answer to Problem 15.1 is negative, then the answer to the following 
problem is positive. 

\noindent{\bf Problem 15.2}
Is there a non SHB link?  

We can define an invariant for (4k+2)-knots 
corresponding to  the $\mu$ invariant for 2-knots. 
We use this invariant and make a similar problem to Problem 15.1.

\footnotesize{
 }

\begin{thebibliography}{900}

\bibitem{Browder} W. Browder: 
The Kervaire invariant of framed manifolds and its generalization, 
{\it Ann. Math.},  90, 157-186, 1969.


\bibitem{CappellShaneson} S.Cappell and J.Shaneson: Link cobordism, 
{\it  Comment.   Math.   Helv.}, 1980.


\bibitem{CO}
T.D.Cochran and K. E. Orr: Not all links are concordant to boundary links 
{\it Ann. Math.}, 138, 519--554, 1993. 

\bibitem{GL}  P. Gilmer and C. Livingston: 
The Casson-Gordon invariant and link concordance, {\it Topology}, 
31, 475-492, 1992.

\bibitem{G} C.McA.Gordon: 
Some higher-dimensional knots with the same homotopy group, 
{\it Quart.J.Math.Oxford } 24, 411-422, 1973. 


\bibitem{HK} O. G. Harrold and S. Kinoshita: 
A theorem on $\theta$-curves and its application 
to a problem of T.B.Rushing, 
{\it Bull.de L'Acad. Polon. des.Sci.}, 28, 631-634, 1980.


\bibitem{Hitt} L.R.Hitt: 
Examples of higher dimensional slice knots which are not ribbon knots,  
{\it Proc. Amer. Math. Soc.}, 77, 291-297, 1979.


\bibitem{Kf} L. Kauffman: On knots, {\it Ann of math studies}, 
 115,  1987.


\bibitem{Kawauchi} A. Kawauchi: 
On links not cobordant to split links, 
{\it  Topology} 19, 321-334,  1980.




\bibitem{Ke} M. Kervaire:  
Les noeudes de dimensions sup\'ereures,  
{\it  Bull.Soc.Math.France}  93, 225-271, 1965.  


\bibitem{KM} M. Kervaire and J. Milnor: Groups of homotopy spheres I, 
{\it Ann. Math. },  77,  504-537, 1963.


\bibitem{K}  S. Kinoshita: 
On $\theta_n$-curves in ${\bf R^3}$ and their constituent knots,  
{\it In:Topology and Computer sciences, Kinokuniya Co.Ltd}, 211-216, 1987.


\bibitem{KT} S. Kinoshita and H. Terasaka: 
 On union of knots, Osaka J. Math, 9, 131-153, 1957. 
 
 
\bibitem{Kirby}
R. Kirby:  The topology of 4-manifolds 
{\it Lecture Notes in Math  (Springer Verlag) }
vol. 1374, 1989

 
 \bibitem{KL}  P. Kirk and C. Livingston:  
Vassiliev invariant of two component links and the Casson-Walker invariant,  
{\it Topology},  36,  1333-1353,  1997


\bibitem{L1}  J. Levine:   
 Polynomial invariant of knots of codimension two, {\it  Ann. Math.}, 
84,  537-554,  1966.

\bibitem{L2}  J. Levine:  Knot cobordism in codimension two, 
{\it Comment. Math. Helv.},  44,  229-244, 1969. 


\bibitem{L3}  J. Levine:  Link invariant via the eta-invariant, 
{\it Comment. Math. Helv.},  
 69, 82-119,  1994. 

\bibitem{Mio}  W, Mio:  On boundary link cobordism, 
{\it  Proc. Camb. Phil. Soc.}, 
101,  259-266,  1989. 


\bibitem{Murasugi}  K. Murasugi:  On the {Arf} invariant of links
{\it  Proc. Camb. Phil. Soc.},  95,  421-434,  1984.



\bibitem{O5}E, Ogasa:
Some results on knots and links in all dimensions. 
{\it Surgery and geometric topology (Edit. A.Ranicki and M. Yamasaki, 
Sakado,1996). Sci. Bull. Josai Univ.} 1997, Special issue no. 2, 81--90. 



\bibitem{O1}  E. Ogasa:  
Nonribbon 2-links all of whose components are trivial knots  
and some of whose band-sums are nonribbon knots
{\it  J. Knot Theory Ramifications}, to appear


\bibitem{O2} E. Ogasa: 
 Ribbon-moves of 2-links preserve the $\mu$-invariant of 2-links, 
{\it math.GT/0004008 in http://xxx.lanl.gov,  UTMS 97-35.}

\bibitem{O3} E. Ogasa:
 Intersectional pair of $n$-knots, local moves of $n$-knots, 
and their associated invariants of $n$-knots
{\it Mathematical Research Letters}, 5, 577-582, 1998, 
 UTMS 95-50 

\bibitem{O4} E. Ogasa:
 Link cobordism and the intersection of slice discs
{\it Bulletin of the London Mathematical Society}, 31,  1-8, 1999. 

\bibitem{O5} E. Ogasa:
 Link cobordism and the intersection of codimension two submanifolds 
{\it In preparation}, which is  
based on a part of the author's Doctor thesis 
(Math. Dept. University of Tokyo 1996). 

\bibitem{Rolfsen} D. Rolfsen: Knots and links
{\it Publish or Perish, Inc.} 1976 


\bibitem{R}  D. Ruberman: 
Doubly slice knots and the Casson-Gordon invariants 
{\it Trans. Amer. Math. Soc.}, 279, 569-588, 1983.


\bibitem{Sato}  N. Sato:  Cobordisms of semi-boundary links,  
{\it  Topology Appl.}, 18, 225-234, 1984.


\bibitem{Si} M. Saito:  On the unoriented Sato-Levine invariant, 
{\it  J. Knot Theory Ramifications}, 2,  335-358, 1993.

\bibitem {Y} A. Yasuhara:  On higher dimensional $\theta$-curves
{\it Kobe J. Math.}, 8, 191-196, 1991

\end{thebibliography}
\end{document}